\documentclass[openbib,twoside,notitlepage]{article}  

\usepackage{amssymb}
\usepackage{latexsym}
\usepackage{exscale}
\usepackage{epsf}
\usepackage{theorem}
\usepackage[dvips]{graphicx}

\newfont{\bsym}{cmbsy10 scaled\magstep2}
\newfont{\bsymi}{cmbsy10}
\newfont{\bmath}{cmmib10 scaled\magstep2}
\newfont{\titfont}{cmbx10 scaled \magstep3}
\newfont{\eighrm}{cmr8}
\newfont{\ack}{cmr10}
\newfont{\Ack}{cmbx10 at 14pt}
\newtheorem{theorem}{Theorem}[section]
\newtheorem{lemma}{Lemma}[section]
\newtheorem{definition}{Definition}[section]
{\theorembodyfont{\rmfamily} 
\newtheorem{remark}{Remark}[section]}

\newcommand{\vsp}{\vspace{12pt}}

\newcommand{\R}{{\mathbb{R}}}

\newcommand{\LL}{{\bf L}}
\newcommand{\ve}{\varepsilon}

\newcommand{\tv}{\hbox{Tot.Var.}}

\newcommand{\n}{\noindent}

\newcommand{\be}{\begin{equation}}
\newcommand{\ee}{\end{equation}}

\newcommand{\lu}{{{\bf L}^1}}

\newcommand{\li}{{{\bf L}^{\infty}}}

\newcommand{\ov}{\overline}

\newcommand{\A}{{\cal A}}

\begin{document}

\setlength{\voffset}{-0.5in}

\title{
{\huge Some Results on }\\
{\huge the Boundary Control of}\\
{\huge Systems of Conservation Laws}\\
      \vspace{0.1in}}

\author{
        {\scshape Fabio Ancona}  \thanks{Dipartimento di Matematica and 
                                 C.I.R.A.M., \
                                 P.zza \ Porta \ S. Donato, \ n.~5, \
                                 40123 - Bologna, \ Italy;
                                 \newline 
                                 E-mail: \texttt{ancona@ciram3.ing.unibo.it.
                                 \vspace{6pt}}}
\and    
        {\scshape Alberto Bressan} \thanks{SISSA-ISAS,
                                via Beirut 2-4, 34014 - Trieste, \ Italy;      
                                      E-mail: \texttt{bressan@sissa.it.}} 
                                \vspace{6pt}
\and    
       {\scshape Giuseppe Maria Coclite} \thanks{SISSA-ISAS,
                               via Beirut 2-4, 34014 - Trieste, \ Italy;      
                                     E-mail: \texttt{coclite@sissa.it.}}
                                \vspace{15pt}}

\date{}

\maketitle
\thispagestyle{empty}

%\vspace{-0.6cm}
%
\begin{abstract}
\vspace{5pt}

\noindent
This note is concerned with the study of 
the initial boundary value problem for systems of
conservation laws from the  point of view of control theory,
where the initial data is fixed and the boundary data
are regarded as control functions.
We first consider the problem of controllability
at a fixed time for genuinely nonlinear
Temple class systems, and present a description
of the set of attainable configurations of
the corresponding solutions 
in terms of suitable Oleinik-type estimates.
We next present a result concerning
the asymptotic stabilization near a constant state for general $n\times
n$ systems. Finally we show with an example that in general  one cannot
achieve exact controllability to a constant state
in finite time.

\end{abstract} \vspace{1.5cm}

%\begin{center}
%Ref. S.I.S.S.A. ??/2002/M (April 2002)
%\end{center}
%
\vspace{1.5cm}

1991\textit{\ Mathematical Subject Classification:} 
%  Conservation Laws:
35L65,
%  PDE in connection with Control Problems:
35B37

\textit{Key Words:} hyperbolic systems, conservation laws,
Temple class systems, 
%  Lipschitz semigroup, 
boundary control,
attainable set.
\vfill
\pagebreak
%
%
%	END TITLE PAGE
%

\setlength{\voffset}{-0in}  
\setlength{\textheight}{0.9\textheight}

\thispagestyle{empty}
\null
\vfill
\pagebreak

\thispagestyle{plain}
\setcounter{page}{1}
\setcounter{equation}{0}

\section[]{Introduction} 
\label{section1}
\indent

Consider an $n\times n$ system of  conservation laws on
a bounded interval
% $[a,b]\subset\R$  
\be u_t+f(u)_x=0\qquad\qquad t\geq 0,~~x\in\,]a,b[\,,\label{1.1}\ee 
with the initial condition
\be u(0,x) = \varphi (x), \qquad a\le x\le b, \label{2}\ee 
%   and the implicit boundary conditions
%   \be\psi_i (u(t,a))= \alpha_i(t), \qquad \psi_i (u(t,b))=
%   \beta_i(t),\qquad t>0,\label{3}\ee
%   we will explain the meaning of this conditions in the following.
and a weak form of the Dirichlet boundary conditions
\be
u_i(t,a)= \alpha_i(t), \qquad u_i(t,b)=
\beta_i(t),\qquad t>0\,
\label{3}\ee
(see \cite{dubfiocco, jofiocco, serre} 
and reference therein for several weak formulations
of~(\ref{3})).

We want to study the effect of boundary conditions  
on the solution of (\ref{1.1}) from the point of
view of control theory. Namely, 
following the same approach
adopted in \cite{AM1} for scalar conservation laws,
we take the initial data $\varphi$ fixed, and, regarding
the measurable maps $\alpha_i, \> \beta_i$ in (\ref{3})
as  {\it control functions}, we want to investigate the property 
of the {\it attainable set} for (\ref{1.1})-(\ref{2}), i.e. 
of the set
\be
\A (T) \doteq \Big\{ u(T,\cdot)~; ~~
u \ {\rm \ is \ a \ sol. \ \ to} \ \ (\ref{1.1})-(\ref{3})\Big\}\,,
\label{att}
\ee
% \subset \, L^1 \big(\, [a,b]; \R^n\big),$$
%
which consists of all profiles that can be attained 
at a fixed time $T>0,$ by 
entropy weak
solutions of (\ref{1.1})-(\ref{3}). 
For the definitions and the basic properties of weak solutions
we refer to \cite{B}.
See also \cite{liu1, sable, A, AC} for results concerning the existence 
and $\lu$ stability of entropy weak  
solutions of the mixed problem
taking values in the space $BV$ of functions with bounded
variation,  and \cite{brgoa1, A-G} for the case
of $\li$ solutions (with possibly unbounded variation)
of Temple class systems.

Throughout, we shall assume:    
\begin{itemize}

\item[$(H_1)$] the map $f:\Omega \longrightarrow \R^n$ is smooth and 
$\Omega \subset \R^n$ is open;

\item[$(H_2)$] the system (\ref{1.1})
is strictly hyperbolic, i.e. the Jacobian matrix $Df(u)$
has $n$ real and  distinct eigenvalues 
$\lambda_1 (u) < ...<\lambda_n (u),\> u\in \Omega$
(with a corresponding basis of eigenvectors
$\{r_1(u), \dots , r_n(u)\}$); 

\item[$(H_3)$] each
characteristic field $r_i$ is linearly degenerate or genuinely nonlinear
in the sense of Lax~\cite{lax}; 

\item[$(H_4)$] there exist $p\in
\{1,...,n\}$ and $c_0 >0 $ such that  
\be \lambda_1 (u) < ...<\lambda_p (u) \le
-c_0 <0< c_0 \le \lambda_{p+1} (u)<...<\lambda_n (u),\quad u\in
\Omega.\label{1.2}\ee   
%
%for each $u\in \Omega.$ 

\end{itemize}

By $(H_4)$, for a solution defined on the strip $t\geq 0,~~x\in\,]a,b[\,$, 
there will be $n-p$ characteristics entering at the boundary 
point $x=a$, and $p$ characteristics entering at $x=b$. 
The initial-boundary value problem is thus well posed 
if we prescribe $n-p$ scalar conditions at $x=a$, and $p$ 
scalar conditions at $x=b$ (see \cite{LY}).  

%\begin{center}
%\leavevmode \epsfxsize=3in \epsfbox{pas2.eps}
%\end{center}
%\begin{center} 
%figure 1
%\end{center}

\begin{definition}
\label{excontr} 
Given $\varphi\in L^1([a,b]),\, v\in \Omega$, and $T>0$, 
we say that the problem
(\ref{1.1})-(\ref{2}) is exact controllable at time $T$ to the state $v$
if and only if there exist 
measurable maps $\alpha_i,\beta_i$ such that the solution of  
(\ref{1.1})-(\ref{3}) satisfies 
$$u(T,\cdot)\equiv v,\qquad {\it a.\> e. \> in}\>\> [a,b].$$
\end{definition}

\begin{definition}
\label{ascontr} Given $\varphi\in L^1([a,b])
,\, v\in \Omega$, we say that the problem
(\ref{1.1})-(\ref{2}) is  asymptotic stabilizable near the state $v$
if and only if there exist measurable maps 
$\alpha_i,\beta_i$ such that the solution of  
(\ref{1.1})-(\ref{3}) satisfies 
$$u(t,\cdot)\longrightarrow v,\qquad {\it in}\>\> L^1([a,b])\qquad {\it as}\>\>
t\longrightarrow +\infty.$$ \end{definition}

In this note we present some recent results 
obtained by the authors~\cite{A-C,BC}
concerning both the problem of exact controllability
and of asymptotic stabilization near a constant state.
We will first consider the case of Temple systems~\cite{temple},
for which it was obtained a characterization of
the corresponding attainable set (\ref{att})
in terms of suitable Oleinik-type estimates,
which is a natural extension of the results
in~\cite{AM1, AM2} concerning scalar conservation laws.
For general nonlinear systems, one cannot expect such an 
exact controllability result.  Indeed,
even if all wave-fronts in the initial data exit
from the interval $[a,b]$ within finite time, they can generate
new waves by interacting among themselves. In turn (figure 3), 
further interactions
can produce a sequence of wave-fronts remaining within the interval
$[a,b]$ for all times $t>0$. Therefore, the effect of the initial data
on the solution $u(T,\cdot)$ may never be completely erased, 
no matter how large we choose the terminal time $T$.
Hence, we will present a result concerning the asymptotic
stabilization of a general system of
conservation laws near a constant state.
Finally, we discuss a counterexample to the
exact controllability concerning a class of $2\times 2$ systems
for which, in general, a constant state $u^*$ cannot be
attained, in a finite time $T$.

An outline of these results 
established in
\cite{A-C,BC} is given in the following sections.
\vsp

\section[]{The Attainable Set for Temple Class Systems} 
\label{section2}
\indent

Our first result is concerned with the problem
of exact controllability of Temple class 
systems~\cite{serre, temple}, which are systems
that satisfy the
following additional assumption.

\begin{itemize}
\item[$(H_5)$] There exists a complete set of Riemann coordinates
$w=(w_1,\ldots,w_n)$
such that each level set $\{u\,;~~w_i(u)={\rm constant}\}$
is an hyperplane. 
\end{itemize}

\n As a consequence, all integral curves of 
the eigenvectors are straight lines
and coincide with the Hugoniot curves.
We shall also assume 

\begin{itemize}
\item[$(H_6)$]as $w$ ranges within the product set
$\Gamma\doteq [w_1^-,w_1^+]\times\cdots\times[w_n^-,w_n^+]$, the corresponding
state $u$ remains inside the domain
$\Omega$ and each characteristic field is genuinely nonlinear.
\end{itemize}

In the case of systems
of this type, the well-posedness theory
for the mixed problem was established in \cite{A-G}
within domains of $\li$ functions
(with possibly unbounded variation).
Here, 
the boundary condition
is formulated in terms of the strong $\lu$ trace
of the solution $u$
%  weak trace of the flux $f(u)$ 
at the 
boundary  and, in the same spirit of~\cite{dubfiocco}, is
based on the definition of a time-dependent
set of {\it admissible
boundary data}, that is related to the notion of Riemann problem. 
Moreover, for sake of uniqueness,
it was introduced in \cite{A-C, A-G}  
a definition of {\it entropy admissible weak solution} to the mixed problem
that includes an entropy admissibility condition of
Oleinik type.

Notice that, for such systems, wave interactions 
can only change the speed of wave fronts, without 
modifying  their amplitudes.
Therefore,  the only restriction to boundary controllability 
is the decay  due
to genuine nonlinearity.
We thus consider a set of maps,
defined in 
terms of the above Riemann coordinates, that
satisfy certain  Oleinik-type conditions
on the decay of positive waves.
\be 
\begin{array}{ll}
\!\!\!\!\!&\!\!\!
K^\rho\!\doteq\!
\left\{ \psi\in{\bf L}^{\infty}([a,b],\,\Gamma)~;\> \,
\begin{array}{ll}
\displaystyle{\frac{w_i(\psi(y))-w_i(\psi(x))}{y-x}}
\leq \displaystyle{\rho\over x-a}\quad 
\left\{
\!\!\!\!\!\!\!\!
\begin{array}{ll}
&\hbox{for \, a.e.} \quad a<x<y<b,
\\
\noalign{\smallskip}
&\hbox{if} \quad i\in \{p+1,...,n\}
\end{array}
\right.
\\
\noalign{\medskip}
\displaystyle{\frac{w_i(\psi(y))-w_i(\psi(x))}{y-x}}
\leq \displaystyle{\rho\over b-y}\quad 
\left\{
\!\!\!\!\!\!\!\!
\begin{array}{ll}
&\hbox{for \, a.e.} \quad a<x<y<b,
\\
\noalign{\smallskip}
&\hbox{if} \quad i\in \{1,...,p\}
\end{array}
\right.
\end{array}
\!\!\!\!\!
\right\}.
\\
\noalign{\medskip}
\end{array}
\label{1.6}   
\ee

\n
The inequalities in (\ref{1.6}) 
reflect the fact that positive waves
entering through the boundary at $a$ or at $b$ decay in time.
Therefore, their density is inversely proportional to their distance
from their entrance point on the boundary.

We can now state our first main result (see \cite{A-C}).

\begin{theorem} \label{theorem1}
Let  (\ref{1.1}) be a system of Temple class, and assume
that $(H_1)$- $(H_6)$ are verified.
Then, letting  $\A(T)$ be the attainable set
defined in (\ref{att}) (in which the solution is
understood as an ``entropy admissible weak solution''), 
the following hold:
\begin{itemize}
\item[$(i)$]for every  fixed
$\overline\tau>0,$ there exists $\rho=\rho(\overline\tau)>0$ such that
\be
\A(\tau) \subseteq  K^{\rho},
\qquad\ \  \tau \geq \overline\tau\,;\label{t2}
\ee
\item[$(ii)$]
there exist $T>0$ and $\rho'<\rho(T)$, such that
\begin{eqnarray}
K^{\rho'}  &\subseteq & \A(\tau),
\qquad\, \tau > T\,;
\label{t1}
\end{eqnarray} 
\item[$(iii)$] $\A(T)$ is a compact subset of $L^1 ([a,b])$ for each $T>0$.
\end{itemize}
\end{theorem}

The first property $(i)$ is an immediate consequence
of the definition of entropy admissible solution
that satisfies suitable Oleinik-type estimates.
\vskip 7pt

The proof of $(ii)$ is established in two steps. 

\n {\bf 1) Backward Construction of Front  Tracking Solutions.}
We take 
\be \tau >T\doteq 4\,{b-a\over\lambda^{\min}},
\qquad\quad\lambda^{\min}\doteq \min_i |\lambda_i|\,, \label{tau}\ee 
and, for a given function $\psi \in K^{\rho'}$, we
construct a sequence of approximate solutions $u^\nu$
on the strip $[0, \tau]\times [a,b]$
such that
\be
{\setlength\arraycolsep{2pt}
\begin{array}{ll}
u_\nu(\tau,\cdot)&=\psi_\nu
\qquad\quad
\psi_\nu \ \, \displaystyle{\mathop{-\!\!\!\rightarrow}^{\lu}} \ \, \psi,
\\
u_\nu(0,\cdot)&=\overline u_\nu
\qquad\quad
\overline u_\nu \ \, 
\displaystyle{\mathop{-\!\!\!\rightarrow}^{\lu}} \ \, \varphi,
\end{array}}
\ee
with the following procedure.
We partition the strip $[0, \tau]\times [a,b]$
in three regions. 
On the rectangle $[(3/4)T, \tau]\times [a,b]$,
starting from $t=\tau$, 
we construct backward in time the front tracking solution $u_\nu$
relying on the fact that the Oleinik estimates
of the definition (\ref{1.6}) of~$K^{\rho'}$
guarantee that two rarefaction fronts of the 
same family never cross in $\Omega$  (see figure 1a). 
The total number of wave-fronts in $u^\nu(t, \, \cdot)$
decreases as $t \downarrow (3/4)T$ whenever a (backward) front crosses
the boundary points $x=a, \, x=b$.
Therefore, since 
the maximum time taken by fronts of $u^\nu$ to
cross the interval~$[a,b]$
is $(b-a)/\lambda^{\min}$, the definition~(\ref{tau}) of $T$
guarantees that all the (backward) fronts of $u^\nu$ will hit 
the boundaries $x=a, \, x=b$ within 
some time $\tau'\in\,](3/4)T,\,\tau[$\,.
Hence, there will be some constant state $\omega \in \Omega$
such that $u^\nu((3/4)T, \cdot) \equiv \omega$.
We next define $u_\nu$ on the rectangle 
$[0, T/4]\times [a,b]$ as 
the restriction to  $[0, T/4]\times [a,b]$
of the front tracking solution 
to the Cauchy problem for (\ref{1.1}), with initial data
$$\overline u(x) = \left\{
\!\!\begin{array}{ll}
\ov u^\nu(a+) & \>\textrm{if \ \ $ x<a$},\\ 
\ov u^\nu(x) & \>\textrm{if \ \ $a \le x \le b$},\\
\ov u^\nu(b-) & \>\textrm{if \ \ $ x > b$}.
\end{array} \right.
$$
Since $u^\nu$ contains only fronts originated
at the points of the segment $\{(0,x)\,;\,x\in[a,b]\}$,
because of (\ref{tau}) these wave-fronts 
cross the whole interval $[a, b]$ and exit from the
boundaries $x=a, \, x=b$ before time $T/4$.
Hence, there will be some state $\omega'\in \Omega$
such that $u^\nu(T/4, \cdot) \equiv \omega'$.
Finally, we define $u^\nu(t, \cdot)$ for $t\in [T/4, \, (3/4)T]$
so that $u^\nu(T/4, \cdot) \equiv \omega'$, 
$u^\nu((3/4)T, \cdot) \equiv \omega$.

\parbox{5cm}{
 \begin{center}
\leavevmode \epsfxsize=2.5in \epsfbox{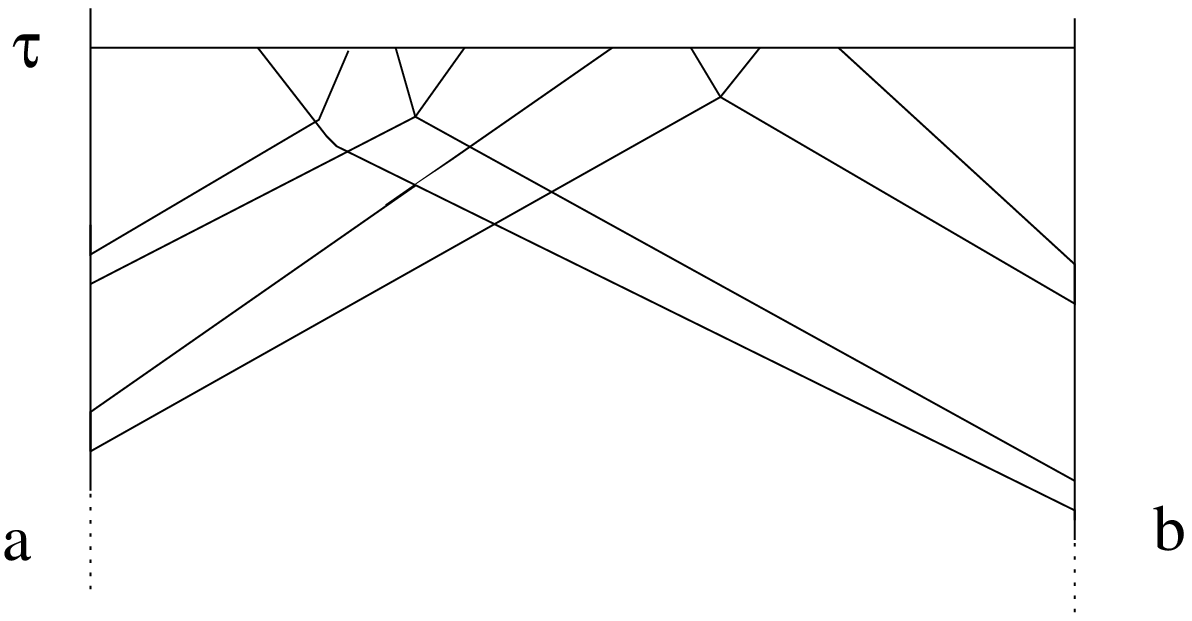}
\end{center}     }
\hfill
\parbox{5cm}{  
 \begin{center}
\leavevmode \epsfxsize=1.7in \epsfbox{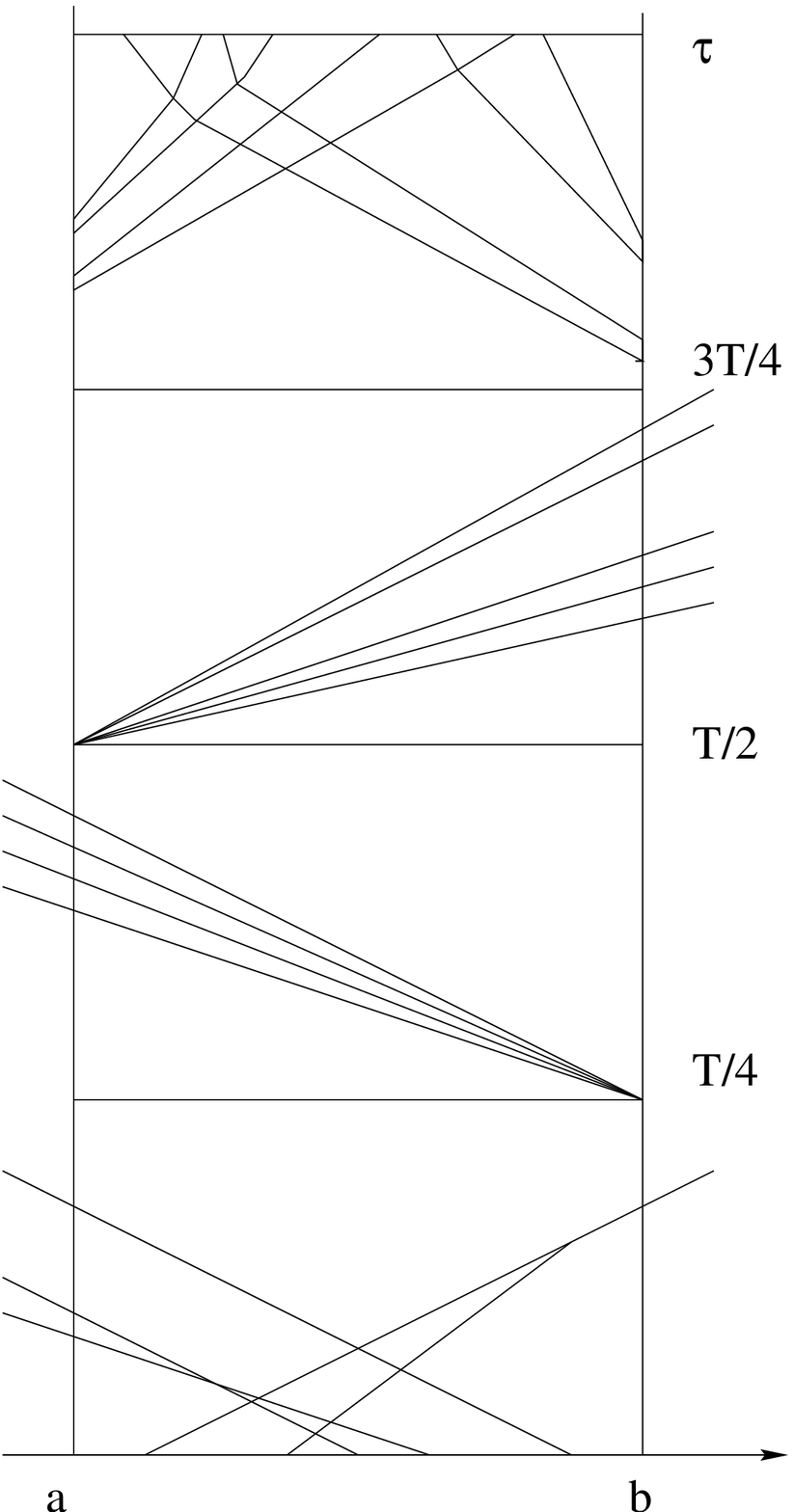}
\end{center}         }
\begin{center} 
~~~~~figure 1a~~~~~~~~~~~~~~~~~~~~~~~~~~~~~~~~~~~~~~~~~~~~~~~~~~~~~~~~~~figure 1b
\end{center}

%\begin{center}
% \vskip0truecm\noindent
% \begin{tabular}[c]{cc}
% \leavevmode 
% \epsfxsize=2.5in \epsfbox{pas5.eps} & \epsfxsize=1.7in
%\epsfbox{pas3.eps}\\ \end{tabular}\\
%\end{center}
%\begin{center} 
%~~~~~~~~~~~figure 1a~~~~~~~~~~~~~~~~~~~~~~~~~~~~~~~~~~~~~~figure 1b
%\end{center}

\n {\bf 2) Convergence of Front  Tracking Solutions.}  
Since the sequence of approximate solutions
constructed by our algorithm 
satisfy the Oleinik-type estimates on the decay of positive waves,
they 
have a
uniformly bounded variation on each interval of
the type $[a+\ve, b-\ve], \> \ve >0$.
Thus, applying Helly's
Theorem, and by  a diagonal procedure, we can extract a
subsequence that converges in $\lu$ to a weak solution $u$
of (\ref{1.1}). 
We then extend a regularity property 
established by A.Bressan and P.G. LeFloch~\cite{brfiocco}
for solutions with small total variation of general
genuinely nonlinear systems, to the case of solutions with
arbitrarily large variation of genuinely nonlinear Temple class
systems.
This property guarantees that, for Temple systems,
solutions of the mixed problem (\ref{1.1})-(\ref{3})
are continuous outsided a countable number of
Lipschitz curves. As an immediate consequence we deduce that
the solution $u$ admits a strong $\lu$ trace
at the boundaries $x=a, x=b$ and satisfies
the corresponding boundary conditions.
\vskip 7pt

Concerning (iii), the compactness of $\A(T)$
is achieved with the same type of arguments used
to establish the convergence of the approximate solutions
in the proof of (ii).

%\begin{center}
%\leavevmode \epsfxsize=1.7in \epsfbox{pas3.eps}
%\end{center}
%\begin{center} 
%figure 3
%\end{center}

%\begin{center}
%\leavevmode \epsfxsize=2in \epsfbox{pas5.eps}
%\end{center}
%\begin{center} 
%figure 2
%\end{center}

%\begin{center}
% \vskip0truecm\noindent
 %\begin{tabular}{cc}
 %\leavevmode 
% \epsfxsize=2.5in \epsfbox{pas5.eps} & \epsfxsize=1.7in
%\epsfbox{pas3.eps}\\ \end{tabular}\\
%\end{center}
%\begin{center} 
%~~~~~~~~~figure 1a~~~~~~~~~~~~~~~~~~~~~~~~~~~~~~~~~~~~~~figure 1b
%\end{center}

\section[]{Asymptotic Stabilization} 
\label{section3}
\indent

We now consider a general $n\times n$ system and show that, starting
with an initial data with small oscillations,
the system can be asymptotically steered to any constant state (see
\cite[Theorem~1]{BC}).

\begin{theorem} \label{theorem2}
Let $K$ be a compact, connected  
subset 
of the open domain 
$\Omega\subset\,\R^n$ and assume that (\ref{1.1}) satisfies $(H_1)-(H_4)$.   
Then there exist constants $C_0,\delta,\kappa>0$ such that the following  
holds. 
For every constant state  
$u^*\in K$ and every initial data $\varphi:[a,b]\mapsto K$ with 
$\tv\{\varphi\}<\delta$, 
there exists an entropy weak 
solution $\> u = u (t, x) \>$ of (\ref{1.1})-(\ref{2}) such that, for all
$t>0$, \be\tv\big\{u(t,\cdot)\big\}\leq C_0\,e^{-2^{\kappa t}}\,,\qquad\big\Vert
u(t,\cdot)-u^*\big\Vert_{L^\infty} \,\leq C_0\,e^{-2^{\kappa  t}}.\label{1.8}\ee 
\end{theorem}

The idea of the proof is as follows. Call $\lambda_*>0$ a lower bound for
the absolute value of all wave speeds and set $\tau\doteq(b-a)/\lambda_*$.
In this way, all waves present in the solution at a given time $t$ will
exit through one the boundaries within time $t+\tau$.
On the first interval $[0,\tau]$ we let all waves exit, arranging
the boundary values at $x=a$ and at $x=b$ so that no reflected
waves ever enter the domain $[a,b]$.  

%\begin{center}
%\leavevmode \epsfxsize=2in \epsfbox{cat3.eps}
%\end{center}
%\begin{center} 
%figure 3
%\end{center}
Therefore, the only waves present in the solution at time $\tau$ are those
generated by interactions, in the interior of the interval $]a,b[$.
They can be estimated as
$$\tv\big\{u(\tau,\cdot)\big\}\leq C_0 \tv\big\{u(0,\cdot)\big\}^2.$$
The above estimate shows that the solution $u(\tau,\cdot)$
remains very close to some constant state, say $u^\dagger$.  
If $u^\dagger\not= u^*$, we 
suitably change the boundary conditions, producing new incoming waves at
$(\tau, b),\>(2\tau, a)$, achieving the bounds 
$$\tv\big\{u(3\tau,\cdot)\big\}\leq C_0 \tv\big\{u(0,\cdot)\big\}^2,\quad
\big\Vert u(3\tau,\cdot)-u^*\big\Vert_{L^\infty} \,\leq C_0\, \big\Vert
u(0,\cdot)-u^*\big\Vert_{L^\infty}^2.$$
Repeating inductively the same strategy in the time intervals
$[3\tau, 6\tau]$, $[6\tau, 9\tau]\ldots$ , we obtain the result.

\begin{center}
\leavevmode \epsfxsize=2in \epsfbox{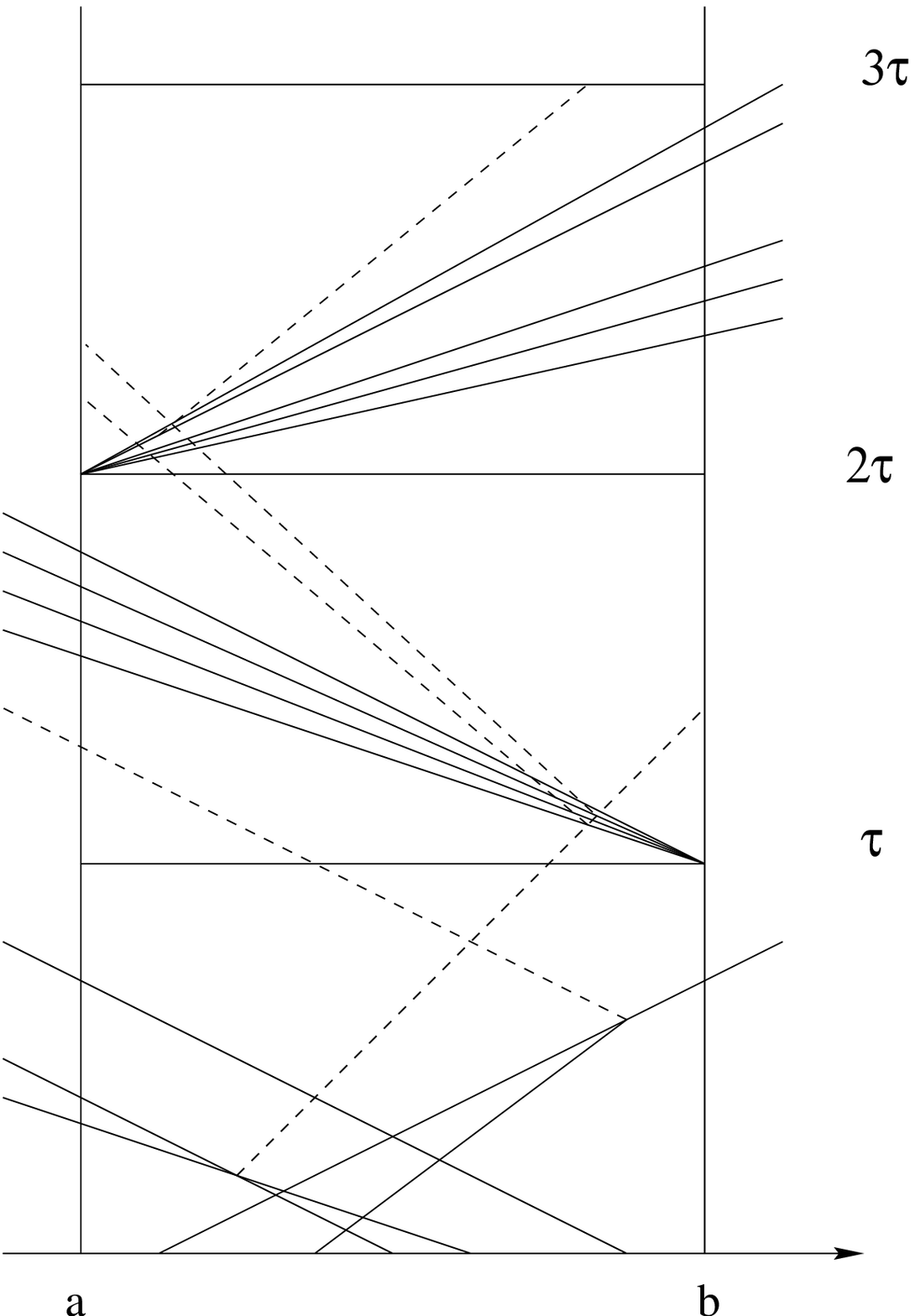}
\end{center}
\begin{center} 
figure 2
\end{center}

\section[]{ A Counterexample to  Exact Controllability} 
\label{section4}
\indent

An interesting question is whether one can reach exactly a constant
state $u^*$ within a finite time $T$.   By Theorem \ref{theorem1}, this is 
certainly the case for Temple class systems.  In the final part of this paper, we show that
this exact controllability cannot be attained in finite time, in general.

Our counterexample is concerned with a class of $2\times 2$ 
strictly hyperbolic, genuinely nonlinear systems
with the property that the interaction of two shocks of the same family
generates a shock in the other family (see \cite[Theorem 2]{BC}).
This is the case for the system (see \cite{DP}):
\begin{equation}
\left\{
\begin{array}{rrl}
\rho_t+(u\rho)_x&=&0,\\[4pt]
u_t+\left( \displaystyle{\frac{u^2}2}+\displaystyle{\frac{K^2}{
\gamma-1}}\rho^{\gamma-1}\right)_{\!\!x}&=&0. \end{array}
\right. \label{1.9}
\end{equation} 
with $1<\gamma<3$. Here $\rho>0$ and $u$ denote the density and the velocity
of a gas,  respectively. 

\begin{theorem} \label{theorem3} Consider a 
$2\times 2$ system of conservation laws satisfying $(H_1)$ and the following

\begin{itemize}
\item[$(H_7)$]  there exist $0 < \lambda_* < \lambda^*$ such that
$$-\lambda^* < \lambda_1(u) <-\lambda_* <0 < \lambda_* < \lambda_2(u)
<\lambda^*,$$
% where $\lambda_1(u) , \lambda_2(u)$ are the
%eigenvalues of  $Df(u)$ and
$$D\lambda_1\cdot r_1 >0,\qquad\qquad
D\lambda_2\cdot r_2>0,$$ 
 $$r_1\wedge r_2<0,\qquad r_1\wedge(Dr_1\cdot
r_1)<0,\qquad  r_2\wedge(Dr_2\cdot r_2)<0,$$ 
  where $r_1(u) , r_2(u)$ are the right eigenvectors of  $Df(u)$.
\end{itemize}

\n Let $\varphi \in BV \big( [a,b]; \R^2\big)$ with small total variation
and a dense set of shocks. Every entropic solution of (\ref{1.1})-(\ref{2})
has a dense set of shocks in $u(t, \cdot),$ for each $t\ge0.$ In particular,
$u(t,\cdot)$ cannot be a constant.
\end{theorem}

\begin{center}
\leavevmode \epsfxsize=2in \epsfbox{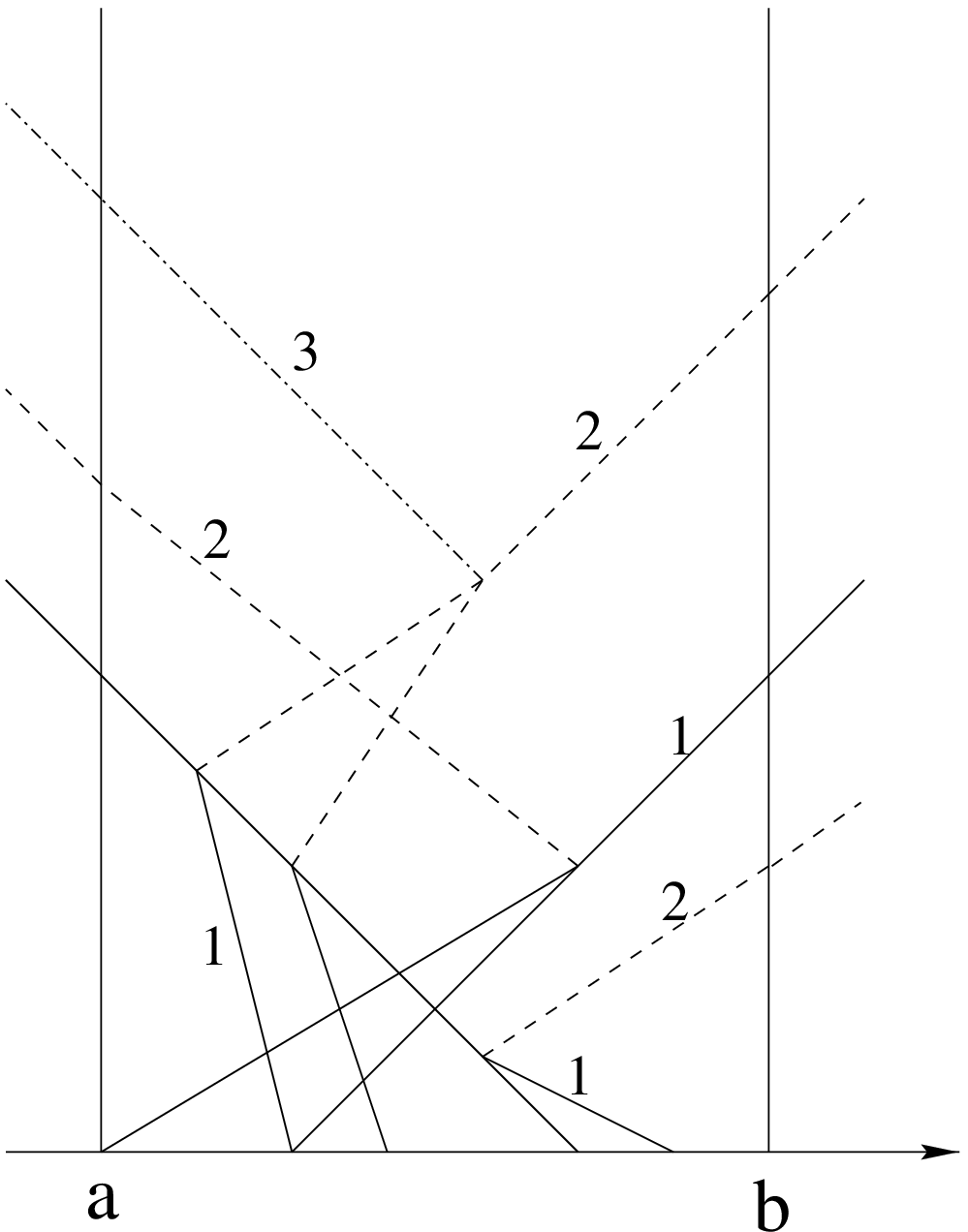}
\end{center}
\begin{center} 
figure 3
\end{center}

\n Toward a proof, we note that the geometric properties 
of the system imply that
\begin{itemize}
\item[a)] the interaction of two shocks of the same family 
produces an outgoing shock of the other family,

\item[b)] the interaction of a shock with a rarefaction wave of the  same
family  produces  a rarefaction wave in the other family. 
\end{itemize}
As in  the scalar case  we have (see
\cite[Section 3]{BC}).

\begin{lemma} 
\label{lemma1}
{\bf (Oleinik type estimate)}  Let $ u= u (t,x)$ be a
solution of (\ref{1.1}) with $ n=2 $ and satisfying  $(H_1)$ and $(H_7)$. There
exist   $k, \delta >0 $ such that, if $$\tv (u (t, \cdot ))< \delta,$$  
 then
$$ \omega_i(t,y) -\omega_i(t,x)\leq { k \over t}\cdot (y-x),\qquad x<y, \>\> 
 t>0,\>\> i =1,2,$$
where $\omega_1,\> \omega_2$ are the Riemann coordinates associated to (\ref{1.1}).
\end{lemma}

By the previous properties, a shock can never be completely canceled by
interactions with rarefaction waves of the same family.  Hence, it
can only disappear by exiting from one of the boundaries $x=a$ or $x=b$.
However, if the set of shocks at time $t=0$ is everywhere dense,
these shocks will interact among each other on a dense set of points
in the domain $[a,b]\times [0,\infty[\,$, and give rise to a dense 
set of new shocks. One can arrange so that the total strength of these shocks
quickly approaches zero, according to Theorem 3.1, but cannot
become exactly zero within finite time.
For all details we refer to 
\cite[Section 3]{BC}.

\begin{remark}
\label{remark1} The previous analysis breaks down in the case of the
$p-$system, because in this case the interaction of 
two shocks of the same family 
produces a centered rarefaction wave of the other family. 
In particular, an
Oleinik type estimate cannot holds. 
\end{remark}

\vsp
\newcommand{\auth}{\textsc}
\newcommand{\tit}{\textrm}
\newcommand{\jou}{\textit}

\end{document}